\documentclass[12pt]{amsart}
\usepackage{amssymb,latexsym,eufrak,amsmath,amscd, graphicx}

\usepackage{amsfonts}
\usepackage{amscd}

\def\ZZ{{\mathbb Z}}
\def\NN{{\mathbb N}}
\def\FF{{\mathbb F}}
\def\CC{{\mathbb C}}

\def\RR{{\mathbb R}}
\def\QQ{{\mathbb Q}}

\def\cI{\mathcal{I}}

\def\ba{\mathbf{a}}

\def\fa{\mathfrak{a}}

\DeclareMathOperator{\lct}{lc}

\DeclareMathOperator{\pt}{c}
\newtheorem{lemma}{Lemma}[section]
\newtheorem{theorem}[lemma]{Theorem}
\newtheorem{corollary}[lemma]{Corollary}
\newtheorem{proposition}[lemma]{Proposition}

\theoremstyle{definition}

\newtheorem{remark}[lemma]{Remark}

\newtheorem{problem}[lemma]{Problem}

\newtheorem{example}[lemma]{Example}

\theoremstyle{remark}
\newtheorem*{remark*}{Remark}
\newtheorem*{note*}{Note}

\begin{document}

\title{Roots of Bernstein-Sato polynomials for monomial ideals:
a positive characteristic approach}

\author[N. Budur]{Nero Budur}
\address{Department of Mathematics, The Johns Hopkins
University, Baltimore, MD 21218, USA}
\email{{\tt budur@math.jhu.edu}}

\author[M. Musta\c{t}\v{a}]{Mircea~Musta\c{t}\v{a}}
\address{Department of Mathematics, University of Michigan,
Ann Arbor, MI 48109, USA}
\email{{\tt mmustata@umich.edu}}

\author[M. Saito]{Morihiko Saito} 
\address{RIMS Kyoto University, Kyoto 606-8502, Japan}
\email{{\tt msaito@kurims.kyoto-u.ac.jp}}

\begin{abstract}
We describe the roots of the Bernstein-Sato polynomial 
of a monomial ideal using reduction mod $p$ and invariants of
singularities in positive characteristic. 
We give in this setting
a positive answer to a problem from \cite{MTW} concerning
the dependence on the characteristic for these invariants of singularities.
\end{abstract}

\thanks{2000\,\emph{Mathematics Subject Classification}.
Primary 14B05; Secondary 13N10, 32S40}
\thanks{The research of the second author was partially supported
by NSF grant 0500127}

\maketitle

\markboth{N. BUDUR, M. MUSTA\c T\v A and M. SAITO}
{ROOTS OF BERNSTEIN-SATO POLYNOMIALS}

\section{Introduction}

The Bernstein-Sato polynomial (or $b$-function) of an
arbitrary ideal in a polynomial ring
was introduced in \cite{BMS} generalizing the case of principal
ideals.
For monomial ideals, it was shown there that this
polynomial can be computed algorithmically in principle,
using for example
\emph{Macaulay}2.
In this paper we use a positive characteristic approach
to give a description of the roots of the Bernstein-Sato
polynomial in this case. This is a by-product of positive answers
to questions in \cite{MTW} on the dependence on the characteristic for some
invariants of singularities.

In order to explain our approach, we recall the definition of
the invariants from \cite{MTW}.
Let $\fa$ be a nonzero ideal in a regular local 
ring $R$ of characteristic $p>0$.  
Suppose that $J$ is
a proper ideal in $R$ whose radical contains $\fa$.
For every positive integer $e$, let $J^{[p^e]}$ be the ideal generated by the
$p^e$th powers of the elements of $J$, and 
$$\nu^J_{\fa}(p^e):=\max\{\ell\geq 0\mid\fa^{\ell}\not\subseteq J^{[p^e]}\}.$$

Suppose now that 
we start with ideals $\fa$ and $J$ in $\ZZ[X_1,\ldots,X_n]$ 
such that $\fa$ is contained in the radical of $J$,
and let $b_{\fa}(s)$ be the Bernstein-Sato polynomial of $\fa$.
We consider the invariants associated to the reductions mod $p$ of
$\fa$ and $J$ around the origin.
It was shown in \cite{MTW} 
that if $p\gg 0$, then $b_{\fa}(\nu^J_{\fa}(p^e))\equiv 0$
(mod $p$) for all $e$ (see also Proposition~\ref{prop1} below).

One expects that
in many cases there are polynomial formulas for the invariants mod $p$,
formulas depending on a suitable congruence of $p$. More precisely,
in good situations there should be a positive integer $N$, and for every $i$
relatively prime to $N$ there should be a polynomial $P_i\in\QQ[t]$
such that $\nu_{\fa}^J(p^e)=P_i(p)$ whenever $p$ is large enough
and $p\equiv i$ (mod $N$).
When this holds,
Dirichlet's
Theorem on the distribution of prime numbers implies that we get
roots of the Bernstein-Sato polynomial over $\QQ$: in fact, each
$P_i(0)$ is such a root
(see Remark~\ref{remark0}).
Our main results are that such formulas can be given for monomial ideals (see
Theorem~\ref{thm1}), and  furthermore that all the roots can be
obtained in this way (see Theorem~\ref{all_roots}). 

The proofs of the above results give in particular a description
of all roots of the Bernstein-Sato polynomial of a monomial idea.
As a consequence we deduce the following description mod ${\mathbb Z}$ of the
roots. Suppose that $\fa$ is a proper nonzero
ideal in $\ZZ[X_1,\ldots,X_n]$ generated by monomials. The Newton
polyhedron $P_{\fa}$ of $\fa$ is the convex hull in $\RR_+^n$ of those
$u$ in $\NN^n$ such that the monomial $X^u$ is in $\fa$. For every facet 
$Q$ of $P_{\fa}$ that is not contained in a coordinate hyperplene,
there is a unique linear function $L_Q$ on $\RR^n$ having
rational coefficients such that $Q=P_{\fa}\cap L_Q^{-1}(1)$.
(Here a facet means a maximal-dimensional face.)
We denote by
$m_Q$ the smallest positive integer such that $m_QL_Q$ has integer coefficients.

\begin{corollary}\label{mod_Z}
The set consisting of the images in $\QQ/\ZZ$ of the roots of the
Bernstein-Sato polynomial of $\fa$ is equal to 
$$\left\{\frac{m}{m_Q}+\ZZ\mid Q\,{\rm facet}\,{\rm of}\,P_{\fa}\,
{\rm and}\, 0\leq m<m_Q
\right\}.$$
\end{corollary}

One can deduce from our results a
description of the roots of the Bernstein-Sato polynomial,
and not just of their classes mod $\ZZ$ (see Remark~\ref{complete_description}). 
Note that
such a description has to be more involved, as the roots
do not depend only on the integral closure of the ideal (or equivalently,
on the Newton polyhedron). 
Another, more explicit
description of the roots of the Bernstein-Sato 
polynomial of monomial ideals (with a direct,
combinatorial proof) will appear in \cite{BMS1}.

\medskip

A few words about the structure of the paper:
in the next section we recall the definition of Bernstein-Sato polynomials
and give some details in the general set-up about the connection between
their roots and the invariants obtained by reduction mod $p$.
In \S 3 we specialize to monomial ideals, and 
give a direct proof of the existence of the 
Bernstein-Sato polynomial in this case. This will be useful later,
as it provides useful information about the roots
(see Remark~\ref{ingred}). The idea is based
on the approach to computing Bernstein-Sato polynomials of monomial ideals
from \cite{BMS}.
In \S 4, we study the invariants of the reduction mod $p$ 
for monomial ideals.
We give a polynomial formula for these invariants
(depending on a suitable congruence of $p$) 
and show that all the roots of the Bernstein-Sato
polynomial can be obtained by our method.
The fourth section is devoted to some examples. In the Appendix we
show that if $\fa$ is a monomial ideal, it is enough to consider only
positive characteristic invariants that correspond to monomial ideals $J$.

\section{Bernstein-Sato polynomials and reduction mod $p$}

We start by recalling some general facts about Bernstein-Sato polynomials.
For proofs and details we refer to \cite{BMS}.

Let $\fa\subseteq (X_1,\ldots,X_n)\subseteq\CC[X]=\CC[X_1,\ldots,X_n]$
be a nonzero ideal,
and let $f_1,\ldots,f_r$ be nonzero generators of $\fa$.
The Bernstein-Sato polynomial $b_{\fa}(s)\in\CC[s]$
of $\fa$ is the monic generator of the ideal consisting of those $b(s)$ 
for which we have a relation 
\begin{equation}\label{b_poly}
b(s_1+\ldots+s_r)\prod_jf_j^{s_j}
=\sum_cP_c\cdot\prod_{j,c_j<0}{s_j\choose -c_j}
\prod_jf_j^{s_j+c_j},
\end{equation}
where the above sum is over finitely many $c\in\ZZ^r$ such that 
$\sum_jc_j=1$, and where $P_c\in \CC[X,\partial_X,s]$ for all $c$.
(Here $\prod_{j,c_j<0}$ means that the product is over the $j$
such that $c_j<0$.)
As usual, if $m>0$, then the notation 
${s_j\choose m}$ stands for $\frac{1}{m!}s_j(s_j-1)\ldots(s_j-m+1)$.

In (\ref{b_poly}) one has to interpret the equality formally. 
Note that if $r=1$, then this equation can be rewritten as
$$b(s)f^s=P(X,\partial_X,s)\cdot f^{s+1},$$
so we recover the usual definition of the Bernstein-Sato polynomial
of a principal ideal (see \cite{Bj}, \cite{kashiwara}). 

The fact that there is a nonzero $b(s)$ as above is proved in \cite{BMS},
where one also shows that it does
not depend on the choice of generators
and that all its roots are negative rational numbers. 
We refer to \cite{BMS} also for
the motivation for the defining formula (\ref{b_poly})
in terms of $V$-filtrations.

In our case the ideal is defined over $\ZZ$, and the $P_c$ can be defined
over $\QQ$ by the following

\begin{proposition}\label{prop0}
If the ideal $\fa$ is defined over a subfield $K$ of $\CC$,
then the $P_c$ in {\rm (\ref{b_poly})} can be defined also over $K$.
\end{proposition}

\begin{proof}
There is a finitely generated $K$-subalgebra $A$ of $\CC$ such that
the $P_c$ are defined over $A$. Choosing a maximal ideal of $A$
and taking the image in the corresponding residue field $K'$ of $A$,
we get a relation in which the $P_c$ are defined over a finite
extension $K'$ of $K$.
We may assume that $K'/K$ is Galois by enlarging
$K'$ if necessary.
Averaging by the action of
the Galois group (where the Galois group
acts on the coefficients of $P_c$),
we get the assertion.
\end{proof}

{} From now on, we consider only ideals in $\ZZ[X]$,
as this will be enough for our purpose.
We will always work in a neighborhood of the origin. Note that there
is a local notion of Bernstein-Sato polynomial, where we require 
(\ref{b_poly}) to hold only in some open subset containing the origin. 
In the rest of this section one could replace $b_{\fa}$ by this local version.
However, in what follows we are interested in monomial ideals
(which are homogeneous), so in this case there will be no distinction
between local and global Bernstein-Sato polynomials.

\smallskip

If $p$ is a prime number, the ideal $\fa$ defines by reduction mod $p$
(and localization)
an ideal $\fa_p$ in $R_p:=\FF_p[X]_{(X_1,\ldots,X_n)}$, where $\FF_p=
\ZZ/p\ZZ$ as usual.
Recall that if $I$
is an ideal in $R_p$ and if $e\geq 1$, then
$I^{[p^e]}=(g^{p^e}\vert g\in I)$.
To simplify the notation, for $\fa\subseteq {\rm Rad}(J)
\subseteq (X_1,\ldots,X_n)$, we will denote by $\nu^J_{\fa}(p^e)$
the largest $\ell$ such that $\fa_p^{\ell}\not\subseteq J^{[p^e]}_p$.

By Proposition~\ref{prop0}, $b_{\fa}$ and all $P_c$ in (\ref{b_poly}) have
coefficients in $\QQ$, and hence in $\ZZ[m^{-1}]$ for some integer $m$.
Moreover, we may assume that for all $c$ that appear in (\ref{b_poly})
and for all $j$ such that $c_j<0$, $(-c_j)!$ divides $m$.
If $p$ does not divide $m$,
then (\ref{b_poly})
will hold also after reduction mod $p$. We will apply 
this equality by letting 
$s_1,\ldots,s_r$ to be nonnegative integers. Note that in this case,
if $c_j<0$ and $s_j+c_j<0$, then ${s_j\choose -c_j}=0$,
so the corresponding term in (\ref{b_poly}) vanishes.

The key to our approach is the following elementary
observation from \cite{MTW}. The proof below uses
Proposition~{\ref{prop0}}, and may be slightly easier than the original one,
although both arguments are essentially the same.

\begin{proposition}\label{prop1}
Let $\fa$ and $J$ be nonzero ideals in $\ZZ[X]$
such that $\fa\subseteq {\rm Rad}(J)\subseteq(X_1,\ldots,X_n)$.
Let $m$ be as above.
If
$p$ does not divide $m$ and $e\geq 1$, then 
\begin{equation}
b_{\fa}(\nu_{\fa}^J(p^e))=0\,\,\text{in}\,\,\FF_p.
\end{equation}
\end{proposition}

\begin{proof}
We can find nonnegative integers $a_1,\ldots,a_r$ such that
$\sum_ja_j=\nu^J_{\fa}(p^e)$ 
and $\prod_jf_j^{a_j}\not\in J_p^{[p^e]}$. On the other hand, the
hypothesis implies $\prod_jf_j^{b_j}\in J_p^{[p^e]}$ if $\sum_jb_j=\nu^J_{\fa}(p^e)
+1$ and $b_j\geq 0$ for all $j$. 

Using (\ref{b_poly}) mod $p$
for $s_j=a_j$ (together with the remark before this proposition)
and the fact that $J_p^{[p^e]}$ 
is an $\FF_p[X,\partial_X]$-submodule of $\FF_p[X]$, we deduce that
$b_{\fa}(\nu^J_{\fa}(p^e))=0$ in $\FF_p$.
\end{proof}

This can be used to give roots of the Bernstein-Sato polynomial
whenever one can solve the following problem from \cite{MTW}.

\begin{problem}\label{prob1}
Under good conditions on
$\fa$ and $J$ there should exist a positive integer $N$ together with
polynomials $P_i$ of degree $e$ in $\QQ[t]$ for 
every 
$i\in (\ZZ/N\ZZ)^\times$ such that
the following holds.
If a prime $p$ is sufficiently large and
$p\equiv i$ (mod $N$), then $\nu^J_{\fa}(p^e)=P_i(p)$.
Moreover, one should be able to choose $N$ to depend only on $\fa$.
\end{problem}

\begin{remark}\label{remark0}
Note that if we have $P_i$ as in the above problem,
then $b_{\fa}(P_i(p))=0$ in $\FF_p$
for infinitely many primes $p$ by Proposition~{\ref{prop1}}
and Dirichlet's Theorem. Therefore $P_i(0)$ is a root of $b_{\fa}$.
Here we may assume that $P_i$ is defined over $\ZZ[m^{-1}]$, replacing
$m$ by a multiple if necessary.
A basic question is which roots can be obtained in this way
(Example~4.1 in \cite{MTW} shows that there
might be roots which are not detected by the above method).
\end{remark}

\smallskip

We mention that the asymptotic behavior of $\nu_{\fa}^J(p^e)$ for
$e\to\infty$ is measured by the $F$-threshold $\pt^J(\fa_p)$ introduced
in \cite{MTW}:
this is defined by
\begin{equation}
\pt^{J}(\fa_p):=\lim_{e\to\infty}\frac{\nu^J_{\fa}(p^e)}{p^e}
=\sup_e\frac{\nu^J_{\fa}(p^e)}{p^e}.
\end{equation}
The set of such numbers for various $J$ form the
jumping coefficients for the test ideals $\{\tau(\fa_p^{\alpha})\}_{\alpha}$
introduced by Hara and Yoshida in \cite{HY}.
On the other hand, in characteristic
zero we have the multiplier ideals of $\fa$, and the corresponding
jumping coefficients (see \cite{lazarsfeld} for the theory of multiplier ideals).
We refer to \cite{MTW} for an overview of results and open questions
relating the jumping coefficients of the multiplier
ideals and the $F$-thresholds for the reduction mod $p$, when $p\gg 0$.

\begin{remark}
As follows from the above discussion, the ``top degree'' part in $\nu^J_{\fa}(p^e)$
is related to the jumping coefficients of the multiplier ideals of $\fa$,
while the ``free term'' is related to the roots of the Bernstein-Sato polynomial.
Recall that there is also a direct connection between these roots
and the jumping coefficients. More precisely, the largest root of $b_{\fa}$
is $-\lct(\fa)$, where $\lct(\fa)$ is the log canonical threshold of $\fa$
(this is the first nonzero jumping coefficient). Moreover, if 
$\lambda\in [\lct(\fa),\lct(\fa)+1)$ is a jumping coefficient, 
then $-\lambda$ is a root of $b_{\fa}$.
These results are
 proved in \cite{kollar} and \cite{ELSV} in the codimension one case,
and in \cite{BMS} in general.
\end{remark}

We turn now to the monomial case, in which we will give an explicit 
description of the whole picture. As we will see, understanding the $F$-thresholds
is quite easy. 
For example, this follows from the result in \cite{HY}
saying that for monomial ideals the test ideals are the same as the
multiplier ideals.
On the other hand, we will see that
the invariants $\nu^J_{\fa}(p^e)$ give much more information.
In particular, by the method described above we will be able to recover all the
roots of the Bernstein-Sato polynomial.

\section{Bernstein-Sato polynomials of monomial ideals}

{}From now on we assume that $\fa$ is an ideal generated by monomials,
so we may take $f_j=X^{a_j}=\prod_iX_i^{a_{i,j}}$ for all $j$,
where $a_j=(a_{1,j},\ldots,a_{n,j})\in\NN^n$. We 
consider the linear forms  $\ell_i(s)=\sum_ja_{i,j}s_j$ on $\ZZ^r$. 

In this section 
we give a direct proof of the existence of the Bernstein-Sato
polynomial in this case. We first make use of the homogeneity of the ideal
to reinterpret equation~(\ref{b_poly}) as in \cite{BMS}. 
We consider, more generally, the ideal $I_\fa\subseteq\CC[s_1,\ldots,s_r]$
consisting of those polynomials $F$ such that
there is an equality of the form
\begin{equation}\label{b_poly2}
F(s_1,\ldots,s_r)f_1^{s_1}\ldots f_r^{s_r}=\sum_cP_c\cdot\prod_{j,c_j<0}
{{s_j}\choose {-c_j}}f_1^{s_1+c_1}\ldots f_r^{s_r+c_r},
\end{equation}
with $c$ and $P_c$ as in (\ref{b_poly}).

On $\CC[X,X^{-1},s]$ 
we consider the $\ZZ^n$-grading given by $\deg(s_j)=0$ and
$\deg(X_i)=e_i$, the $i$th element of the standard basis of $\ZZ^n$.
This induces a grading on $\CC[X,X^{-1},s]\prod_jf_j^{s_j}$
such that $\deg(h\prod_jf_j^{s_j})=\deg(h)$.
On $\CC[X,\partial_X,s]$ we have the corresponding grading with
$\deg(\partial_{X_i})=-\deg(X_i)$ and the action on 
$\CC[X,X^{-1},s]\prod_jf_j^{s_j}$ 
is compatible with the gradings. 

Note that $\prod_jf_j^{s_j}=\prod_iX_i^{\ell_i(s)}$
and $\prod_jf_j^{s_j+c_j}=\prod_iX_i^{\ell_i(s+c)}$.
It follows that in (\ref{b_poly2}) we may assume that
$\deg(P_c)=-(\ell_1(c),\ldots,\ell_n(c))$ for all $c$.
If $P\in\CC[X,\partial_X,s]$ has degree $m=(m_i)$, then
$$P\in\CC[X_1\partial_{X_1},\ldots,X_n\partial_{X_n},s]\cdot\prod_i\xi_i^{|m_i|},$$
where $\xi_i=X_i$ if $m_i\geq 0$, and $\xi_i=\partial_{X_i}$ if 
$m_i\leq 0$. 

For every $c\in\ZZ^r$, let 
$$g_c:=\prod_{j,c_j<0}{{s_j}\choose{-c_j}}\cdot\prod_{i,\ell_i(c)>0}
{{\ell_i(s)+\ell_i(c)}\choose{\ell_i(c)}}.$$

\begin{proposition}\label{ideal}
With the above notation, we have
$$I_\fa=\left(g_c\mid c\in\ZZ^r,\sum_jc_j=1\right)
=\left(g_c\mid c\in\ZZ^r,\sum_jc_j\ge 1\right).$$
\end{proposition}

\begin{proof}
We get the first equality by the above argument, 
considering the action of 
$\partial_{X_k}^{\ell_k(c)}$ on $\prod_iX_i^{\ell_i(s+c)}$.
For the second equality we need to show that
if $c\in\ZZ^r$ and $\sum_jc_j\geq 1$, then $g_c\in I_\fa$. 
We get this by induction on $\sum_jc_j$. Indeed,
if $\sum_jc_j\geq 2$ and
 $c_{j_0}>0$, let $c'_j=c_j$ for $j\neq j_0$ and $c'_{j_0}=c_{j_0}-1$.
Since $\ell_i(c)\ge\ell_i(c')$, it follows that $g_{c'}$ divides $g_c$.
By the induction hypothesis, $g_{c'}$ lies in $I_{\fa}$, hence so does $g_c$.
\end{proof}

It follows from the above description of the generators of $I_\fa$ that
for any irreducible component $\Gamma$ of
$V(I_\fa)_{\rm red}$, there 
are subsets $A\subseteq\{1,\ldots,r\}$ and $B\subseteq\{1,\ldots,n\}$ together with
$\alpha_j$ in $\ZZ_{\geq 0}$ for $j$ in $A$ and $\beta_i$
in $\ZZ_{<0}$
for $i$ in $B$ such that 
\begin{equation}\label{gamma}
\Gamma=\{u=(u_j)\in\CC^r\mid u_j=\alpha_j\,{\rm for}\,j\,{\rm in}\, A,
\ell_i(u)=\beta_i\,
{\rm for}\,i\,{\rm in}\, B\}.
\end{equation}
Moreover, after possibly enlarging $A$ and $B$, 
we may assume that if $j$ is not in $A$ and $u_j$
is constant on $\Gamma$, then $u_j$ is not in $\ZZ_{\geq 0}$ for $u$ in $\Gamma$, 
and if $i$ is not in $B$ and
$\ell_i$ is constant on $\Gamma$, then $\ell_i(u)$ is not in $\ZZ_{<0}$ for 
$u$ in $\Gamma$.

\begin{lemma}\label{condition}
With the above notation and assumption, if $c$ in $\ZZ^r$  
is such that $c_j\geq-\alpha_j$ for all $j$ in $A$ and
$\ell_i(c)\leq-\beta_i-1$ for all $i$ in $B$, then 
$\sum_jc_j\leq 0$.
\end{lemma}

\begin{proof}
By Proposition~{\ref{ideal}}, 
if $\sum_jc_j\geq 1$ then $\Gamma\subseteq g_c^{-1}(0)$.
Therefore the assertion follows from the above maximality assumption.
\end{proof}

\begin{remark}
We see that the converse of the above argument is also true:
if $A$, $B$, the $\alpha_j$ and the
$\beta_i$ are such that the assertion in Lemma~{\ref{condition}}
holds and if $\Gamma$ is given by
(\ref{gamma}), then $\Gamma\subseteq V(I_\fa)_{\rm red}$.
\end{remark}

By definition, the Bernstein-Sato polynomial of $\fa$
is the monic generator of $I_\fa\cap\CC[s_1+\ldots+s_r]$. 
The existence of this polynomial follows from the next proposition.

\begin{proposition}\label{existence}
With the above notation, the intersection $I_\fa\cap\CC[s_1+\ldots+s_r]$
is nonzero. 
\end{proposition}

\begin{proof}
We need to show that for every irreducible component $\Gamma$ of
$V(I_\fa)_{\rm red}$, the map $u=(u_j)\in\Gamma\longrightarrow\sum_ju_j$
is constant. This is equivalent to the assertion that if 
$A$ and $B$ are such that 
$\Gamma$ is given by (\ref{gamma}) and Lemma~\ref{condition}
holds, then for every $u$ in $\CC^r$ such that $u_j=0$ for all $j$ in $A$
and $\ell_i(u)=0$ for all $i$ in $B$, we have $\sum_ju_j=0$.
Since the vector space cut out by these equations
is defined over $\QQ$, we may assume that $u$ is in $\QQ^r$ and moreover, that
it is in $\ZZ^r$.

Suppose that $\sum_ju_j$ is nonzero. 
After possibly replacing $u$ by $-u$, we may assume that $\sum_ju_j>0$.
In this case we get a contradiction by
applying Lemma~\ref{condition} for $c=u$.
This completes the proof of the existence of the Bernstein-Sato polynomial
$b_{\fa}$
in the monomial case. 
\end{proof}

\begin{remark}\label{ingred}
It follows from the discussion before Proposition~\ref{ideal}
that both $b_{\fa}$ and the $P_c$
satisfying (\ref{b_poly2}) have rational coefficients.
The above proof shows that all roots of $b_{\fa}$ are rational,
since $\Gamma$ is an affine linear subspace defined by equations
with rational coefficients.
Moreover, we have obtained the following
description of the roots:
we need to consider all $A$, $B$, $\alpha_j$ and $\beta_i$ such that
the assertion in Lemma~\ref{condition} holds. 
If $\Gamma$ is defined by (\ref{gamma}) and if it is nonempty, then
$\sum_ju_j$ is constant on $\Gamma$, and its value gives a root of $b_{\fa}$.
In addition, all the roots arise in this way.
\end{remark}

\section{Roots of Bernstein-Sato polynomials for monomial ideals}

We study now the
invariants from \S 2 in the case
when $\fa$ is a monomial ideal. Proposition~\ref{prop_A1} in the Appendix
shows that for every ideal $J\subseteq (X_,\ldots,X_n)
\subseteq\ZZ[X_1,\ldots,X_n]$ such that
$\fa\subseteq {\rm Rad}(J)$,
there is a monomial ideal 
$\widetilde{J}$ such that if $p\gg 0$, then 
$\nu_{\fa}^J(p^e)=\nu_{\fa}^{\widetilde{J}}(p^e)$
for every $e\geq 1$. Therefore in order to understand the functions
$\nu_{\fa}^J$ for various $J$, it is enough to consider the case when
$J$ is a monomial ideal, too.
Our goal is to give an affirmative answer to
 Problem~\ref{prob1} in this setting and
to show that all roots of $b_{\fa}$ are given by our method. Moreover,
along the way we will get a description of these roots.

We assume that $\fa$ and $J$ are proper nonzero monomial ideals 
in $\ZZ[X_1,\ldots,X_n]$ such that $\fa$ is contained in the radical of $J$.
Since we deal with monomial ideals, we can define 
the function $\nu_{\fa}^J$
without taking the reduction mod $p$.
If $q$ is an arbitrary positive integer
(not necessarily a prime power), we put
$$J^{[q]}:=(X^{qw}\vert X^w\in\ J),$$ 
where $X^w=\prod_i X_i^{w_i}$ for $w=(w_1,\dots,w_n)\in\NN^n$.
We also define
$$\nu^J_{\fa}(q):=\max\{t\geq 0\mid\fa^t\not\subseteq J^{[q]}\}.$$

Our first goal is to prove the following theorem. 
Note that it immediately gives a positive answer to Problem~\ref{prob1}
for monomial ideals.

\begin{theorem}\label{thm1}
If $\fa$ is a nonzero proper monomial
ideal, then there is a positive integer $N$ with the
following property. If $J$ is a monomial ideal whose radical contains $\fa$,
then there are rational numbers $\alpha>0$ and $\gamma_j$ for $j=0,\ldots,N-1$, 
such that $\nu^J_{\fa}(q)=\alpha q+\gamma_j$
if $q\equiv j$ (mod $N$) and $q$ is large enough.
\end{theorem}

\begin{remark}\label{remark1}
With the notation in the theorem, note that if $j$ and $N$ 
are relatively prime, then by
Proposition~\ref{prop1} together with Dirichlet's Theorem,
$b_{\fa}(\gamma_j)\equiv 0$ (mod $p$)
for infinitely many primes $p$. Therefore $\gamma_j$ is a root of $b_{\fa}$
for every such $j$.
Note also that if $p\gg 0$ is a prime, then $\alpha$ is equal 
to the $F$-pure threshold 
$\pt^{J}(\fa_p)$.

The proof of Theorem~\ref{thm1}
will give, in fact, a description of $\alpha$ and of the
$\gamma_j$.
Moreover, Theorem~\ref{all_roots} below will show that all roots of
$b_{\fa}$ are given by some $\gamma_j$ as above, for 
a suitable $J$ and some $j$
relatively prime to $N$.
\end{remark}

\bigskip

We start with some preparations.
Note first 
that we may write $J=\bigcap_{i=1}^dJ_i$, where each $J_i$ is generated
by powers of variables $X_j^{b_j}$ for $j\in I_i\subseteq\{1,\ldots,n\}$.
This can be checked, for example, by induction on $n$.
If Theorem~\ref{thm1} holds for every $J_i$,
then it holds for $J$: it is clear that we have 
$\nu_{\fa}^J(q)=\max_i\nu^{J_i}_{\fa}(q)$ for all $q$.
The assertion follows from the fact that if 
a finite set $S$ is such that for every $i$ in $S$ the function
$m\to h_i(m)$ is affine linear for $m$ large enough
and in a suitable congruence class, then so is
the function $m\longrightarrow\max_{i\in S}h_i(m)$.
Therefore from now on we may assume that $J=(X_i^{b_i}\vert i\in I)$
for some $I\subseteq\{1,\ldots,n\}$.

Recall that we denote the generators of $\fa$ by $X^{a_1},\ldots,X^{a_r}$,
where $a_j=(a_{1,j},\ldots,a_{n,j})$. 
We have  $\fa^{t}\subseteq J^{[q]}$ if and only if for 
all $\beta=(\beta_j)\in\NN^r$ with $\sum_j\beta_j=t$, there is $i$ in $I$
such that
$qb_i\leq\sum_{j=1}^ra_{i,j}\beta_j$.
Hence
$\nu^J_{\fa}(q)=\tau(I;(qb_i-1)_{i\in I})$, where
for $w=(w_i)\in\NN^{|I|}$, we put
$$\tau(I;w):=\max\left\{\sum_{j=1}^r\beta_j\vert\beta\in\NN^r,
\ell_i(\beta)\leq w_i\,{\rm for}\,i\, {\rm in}\, I\right\}.$$
Recall that $\ell_i(\alpha)=\sum_{j=1}^ra_{i,j}\alpha_j$.
We denote by $\ell_i$ also the extension of this linear function to $\QQ^r$.

In order to simplify the notation we show that we may assume
$I=\{1,\ldots,n\}$, after replacing, if necessary,
 $J$ by $J':=J+(X_i^M\vert i\not\in I)$
for $M$ sufficiently large.
Indeed, note first that if $b'$ is in $\NN^{|I|}$,
then $\tau(I;b')$ is finite, so the set
$$\{\beta\in\QQ_+^r\mid\ell_i(\beta)\leq b'_i\,{\rm for}\,i\,{\rm in}\,I\}$$
is bounded.
This implies that if $M\gg 0$, then for every $q\geq 1$ and every
$\beta\in\NN^n$ such that $\ell_i(\beta)< qb_i$ for all $i$ in $I$,
we have $\ell_i(\beta)< qM$ for all $i\not\in I$. For such $M$
we have $\nu^J_{\fa}(q)=\nu^{J'}_{\fa}(q)$ for every $q\geq 1$. 

{}From now on we assume that $I=\{1,\ldots,n\}$, and we
put $\tau(w)$ for $\tau(I;w)$.
For every  $w$ in $\QQ_+^n$ we define also
$$\tau_{\QQ}(w):=\max\left\{\sum_{j=1}^r\alpha_j\vert\alpha\in\QQ_+^r,
\ell_i(\alpha)\leq w_i\,{\rm for}\,{\rm all}\,i\right\}.$$
Here we may replace $\alpha\in\QQ_+^r$ with $\alpha\in\RR_+^r$
because $w$ is in $\QQ_+^n$ and the $\ell_i$ have coefficients in $\QQ$.

Let $P_{\fa}$ be the Newton polyhedron of $\fa$ and let
$\Delta$ be the fan decomposition of $\QQ_+^n$ whose cones are the
closed convex cones over the faces of $P_{\fa}$. A maximal such cone
corresponds to a facet $Q$ of $P_{\fa}$ that is not contained
in a coordinate hyperplane, and we denote
by $L_Q$ or by $L_{\sigma}$ the linear function
such that $Q=P_{\fa}\cap L_{Q}^{-1}(1)$.

\begin{lemma}\label{lemma1}
Let $P_{\fa}$ be the Newton polyhedron of $\fa$. 
If there is $\lambda>0$ such that $w$ is in $\lambda P_{\fa}$, then
$$\tau_{\QQ}(w)=\max\{\lambda>0\mid w\in \lambda P_{\fa}\},$$
and if there is no such $\lambda$, then $\tau_{\QQ}(w)=0$. In particular,
$\tau_{\QQ}$ is piecewise linear on $\QQ_+^n$: if 
$\sigma$ is a maximal cone in $\Delta$, then $\tau_{\QQ}=L_{\sigma}$
on $\sigma$.
Moreover, there is a positive integer $N$ such that $\tau_{\QQ}(w)=
\tau(w)$ if  $w_i/N\in\NN$ for all $i$.
\end{lemma}

\begin{proof}
Since $P_{\fa}={\rm conv}(a_1,\ldots,a_r)+\RR_+^n$,
where ${\rm conv}(a_1,\ldots,a_r)$ denotes the convex hull of
the $a_i$, we see that $\frac{1}{\lambda}w$ is in $P_{\fa}$ if and only if
there are $\beta_1,\ldots,\beta_r\geq 0$ such that
$\sum_j\beta_j=1$ and $w_i/\lambda\geq \ell_i(\beta)$
for all $i$,
where $\beta=(\beta_1,\ldots,\beta_r)$. 
This gives the first assertion.

Let $w$ be in the closed cone $\sigma$ over a facet $\sigma^{(1)}$ of 
$P_{\fa}$. If there is $\lambda>0$ such that $w$ is in $\lambda P_{\fa}$,
then for the largest such $\lambda$ we have $w$ in $\lambda\sigma^{(1)}$.
Therefore $\lambda=L_{\sigma}(w)$. Since both $L_{\sigma}$ and $\tau_{\QQ}$
are continuous, we deduce that $L_{\sigma}=\tau_{\QQ}$ on $\sigma$.

We denote by $e_1,\ldots,e_n$ the standard basis of $\RR^n$.
Every maximal cone $\sigma$ in $\Delta$
is generated as a convex cone by
some of the $a_j$ (the vertices of $\sigma^{(1)}$) and some of the vectors
$e_1,\ldots,e_n$.
By Carth\'{e}odory's
theorem (see Proposition~1.15 in 
\cite{Zi}) $\sigma$ can be written as a union of cones, each generated by
a subset of the generators of $\sigma$ that forms a basis of $\RR^n$.

Suppose now that $w$ lies in the maximal cone $\sigma$.
By the above discussion, we can write
$$w=\sum_{j\in I_1}\alpha_ja_j+\sum_{i\in I_2}\beta_ie_i\quad
{\rm with}\,\alpha_j\,{\rm and}\,\beta_i\,{\rm in}\,\RR_+,$$
where the $a_j$ and the $e_i$ lie in $\sigma$ and give a basis of $\RR^n$.
For the last assertion we need to find $N$ such that
for every $w$ as above lying in $(N\NN)^n$, the
$\alpha_j$ and the $\beta_i$ are integers.

Let $N$ be a positive number which is divisible by the determinant of
any square submatrix of $(a_{i,j})$. 
Using the determinant of the submatrix $(a_{i,j})$, with 
$i$ in the complement of $I_2$ and $j$ in $I_1$, we see
that $\alpha_j$ is an integer for every $j$ in $I_1$.
This in turn implies that all $\beta_i$ are integers, and completes the
proof of the lemma.
\end{proof}

With $N$ as in the above lemma, we have a good
understanding of $\tau$ on $(N\NN)^r$.
We describe now the behavior of $\tau$ 
on the congruence classes modulo the subgroup $(N\ZZ)^n$.
Let $\Delta$ be the
fan decomposition in Lemma~\ref{lemma1}. Consider a cone
$\sigma$ in $\Delta$ which is not contained in any of the coordinate hyperplanes, 
and a translate $v+\sigma$ of this cone, for some
$v\in\ZZ^n$.

\begin{lemma}\label{lemma2}
With the above notation, there is $w$ in $\sigma\cap\NN^n$ 
such that the restriction 
of $\tau$ to the intersection of each congruence class with $v+w+\sigma$
is given by an affine linear function. More precisely, for every 
$c=(c_1,\ldots,c_n)$ in $\{0,\ldots,N-1\}^n$ there is $A_c$
in $\QQ_{+}$ such that
for every $u$ in $(v+w+\sigma)\cap\NN^n$
with $u_i\equiv c_i$ (mod $N$) for all $i$, we have
$\tau(u)=\tau_{\QQ}(u)-A_c$.
\end{lemma}

\begin{proof}
Note that we may replace at any time $v$ by $v+v'$ for some $v'$
in $\sigma\cap\NN^n$.
If $\widetilde{\sigma}$ is a maximal cone in $\Delta$ such that
$\sigma$ is a face of $\widetilde{\sigma}$, then by taking $v'$ to be
a large enough multiple of an element in the interior of $\widetilde{\sigma}$,
we may assume that $v+\sigma$ is contained in $\widetilde{\sigma}$.

It is clear that $\tau$ is concave: if $b$, $b'\in\NN^n$
we have $\tau(b+b')\geq\tau(b)+\tau(b')$. 
On the other hand, $\tau(b)$ and $\tau_{\QQ}(Nb)=N\tau_{\QQ}(b)$
are in $\NN$, so
$$0\leq\tau_{\QQ}(b)-\tau(b)\in\frac{1}{N}\ZZ.$$
It follows that given $c$, 
we may choose $w_c\in\sigma\cap\NN^n$ such that
$$\tau_{\QQ}(v+w_c)-\tau(v+w_c)=$$
$$\min\{\tau_{\QQ}(v+w)-\tau(v+w)\mid
w\in\sigma\cap\NN^n, v_i+w_i\equiv c_i\,({\rm mod}\,N)\,{\rm for}\,{\rm all}\,i\}.$$
If $w'\in \sigma\cap (N\NN)^n$, then by concavity we have
\begin{equation}\label{eq9}
\tau_{\QQ}(v+w_c+w')-\tau(v+w_c+w')\leq \tau_{\QQ}(v+w_c)-\tau(v+w_c)
\end{equation}
(note that $\tau_{\QQ}$ is linear on $\widetilde{\sigma}$).
By minimality, we have equality in (\ref{eq9}). 
If we take $w=\sum_cw_c$, then we can find for every $c$
an $A_c$ as required by the lemma.
\end{proof}

\begin{remark}
Note that if $\sigma$ is a maximal cone, then for every $c$ there are infinitely 
many $u$ in $(v+w+\sigma)\cap\NN^n$ such that $u_i\equiv c_i$ (mod $N$)
for every $i$. This is not necessarily the case if $\sigma$ is not maximal.
However, if given $c$ there is one such $u$, then there are infinitely many
with the same property.
\end{remark}

We can solve now the monomial case of Problem~\ref{prob1}.

\begin{proof}[Proof of Theorem~\ref{thm1}]
We use the notation in Lemmas~\ref{lemma1} and 
\ref{lemma2}. We have seen that we may assume 
$J=(X_i^{b_i}\vert 1\leq i\leq n)$.
Consider $N$ and the fan $\Delta$ 
in Lemma~\ref{lemma1}. Let $\sigma$ be the cone in $\Delta$ such
that $b=(b_i)$ lies in the 
relative interior of $\sigma$ (note that since $b_i>0$ for every $i$,
$\sigma$
is not contained in any coordinate hyperplane).
We put $e=(1,\ldots,1)\in\NN^n$ and  let $\widetilde{\sigma}$ be a maximal cone
in $\Delta$
such that $qb-e\in\widetilde{\sigma}$ for $q\gg 0$.
It follows that $\sigma$ is a face
of $\widetilde{\sigma}$. Recall that we have
a linear function $L_{\widetilde{\sigma}}$ 
 whose restriction to $\widetilde\sigma$
is equal to $\tau_{\QQ}$.

Lemma~\ref{lemma2} implies
that we can find 
$v_{\sigma}$ in $\sigma\cap\NN^n$ and for every $c\in\{0,\ldots,N-1\}^n$, a
nonnegative rational number
$A_c^{\sigma}$
such that 
$$\tau_{\QQ}(u)-\tau(u)=A^{\sigma}_c$$ 
if $u$  is in $(v_{\sigma}-e+\sigma)
\cap\NN^n$
and $u_i\equiv c_i$ (mod $N$) for all $i$.

Recall that $\nu^J_{\fa}(q)=\tau(qb-e)$. Moreover,
if $q$ is large enough, then $qb-e$  lies in $v_{\sigma}-e+\sigma$. Given 
$j\in\{0,\ldots,N-1\}$, we take $c=(c_i)$ such that $jb_i-1\equiv c_i$
(mod $N$) for all $i$. If we put $\alpha=L_{\widetilde{\sigma}}(b)=\tau_{\QQ}(b)$
and $\gamma_j=L_{\widetilde{\sigma}}(-e)-A_c^{\sigma}$, 
then the requirement of the theorem is satisfied.
\end{proof}

\begin{remark}\label{complete_description}
For future reference, we give explicitly
the description of the roots of $b_{\fa}$ that are obtained by
our method (Theorem~\ref{all_roots} below shows that these are, indeed,
all the roots of $b_{\fa}$).
For every cone $\sigma$ in our
fan $\Delta$, such that $\sigma$
  is not contained in a coordinate hyperplane, let us choose
a maximal cone $\widetilde{\sigma}$ in $\Delta$ 
with the property that for some $v$ in $\sigma$
we have $v-e+\sigma\subseteq\widetilde{\sigma}$ (hence $\sigma$ is a face
of $\widetilde{\sigma}$). Let $L_{\widetilde{\sigma}}$ be the linear function whose
restriction to $\widetilde{\sigma}$ is equal to $\tau_{\QQ}$.

We consider now those $c$ in $\{0,\ldots,N-1\}^n$ such that there is 
$b$ in $\NN^n$
in the relative interior of $\sigma$ with $b_i-1\equiv c_i$ (mod $N$) for all $i$
(if $\sigma$ is maximal, then all $c$ satisfy this condition). 
With $A^{\sigma}_c$ as in the proof of Theorem~\ref{thm1}, we deduce from
Remark~\ref{remark1} that $L_{\widetilde{\sigma}}(-e)-A_c^{\sigma}$
is a root of $b_{\fa}$. Indeed, it is enough to consider $j$ with
$j\equiv 1$ (mod $N$) and $J=(X_1^{b_1},\ldots,X_n^{b_n})$.
In addition, every root we obtain by our method 
is of this form.
Note that the roots we obtain for $\sigma\in\Delta$ do not depend
on the choice of $\widetilde{\sigma}$: if $\widetilde{\sigma}'$
is another maximal cone that satisfies the same property, then
$L_{\widetilde{\sigma}}$ and $L_{\widetilde{\sigma}'}$ agree on the
linear span of $\sigma$ and $e$.
\end{remark}

\begin{remark}\label{rem_for_cor}
With the notation in
the previous remark,
 the class of the root $L_{\widetilde{\sigma}}(-e)-A_c^{\sigma}$
in $\QQ/\ZZ$ is equal to the class of $L_{\widetilde{\sigma}}(-e)-
L_{\widetilde{\sigma}}(qw-e)=-qL_{\widetilde{\sigma}}(w)$, where 
$w\in\NN^n$ in the relative interior of $\sigma$ and
$q\gg 0$ are such that $w_i-1\equiv c_i$ (mod $N$) for all $i$ and
$q\equiv 1$ (mod $N$). Recall that $N L_{\widetilde{\sigma}}(w)$
is an integer for all such $w$. We see that in order
to compute the set of all such classes, when $\sigma$ and $c$ vary,
it is enough to consider only the maximal cones $\sigma$.
The set we get in this way is the set of classes of
$$\{-\tau_{\QQ}(w) \mid w\in (\ZZ_{>0})^n\}.$$
\end{remark}

\begin{remark}
It follows from Remark~\ref{complete_description}
that
we get the same roots of the Bernstein-Sato polynomial
if we consider only the invariants
$\nu^J_{\fa}(p)$ for ideals $J$ of the form
$(X_1^{b_1},\ldots,X_n^{b_n})$, and 
for $p$ prime and large enough with $p\equiv 1$ (mod $N$).
\end{remark}

The following theorem shows that all the roots of
the Bernstein-Sato polynomial of a monomial ideal
are detected by our method.

\begin{theorem}\label{all_roots}
For every nonzero monomial ideal $\fa$, and for every root
$\lambda$ of $b_{\fa}$ there is a monomial ideal $J$ together with
a rational number $\alpha$ and a positive integer $N'$
such that $\nu^J_{\fa}(q)=\alpha q+\lambda$ for $q$ sufficiently large
and with $q\equiv 1$ (mod $N'$).
\end{theorem}

\begin{corollary}
The procedure described in Remark~\ref{complete_description}
gives all the roots of the Bernstein-Sato polynomial $b_{\fa}$. 
\end{corollary}

\begin{proof}
If $\lambda$ is a root of $b_{\fa}$, let $J$, $\alpha$ and $N'$ be as in
Theorem~\ref{all_roots}. Applying Theorem~\ref{thm1} for $\fa$ and $J$,
we see that there are $\alpha'$ and $\lambda'$ such that
$\nu_{\fa}^J(q)=\alpha' q+\lambda'$ if $q$ is large enough and
$q\equiv 1$ (mod $N$). Taking $q\equiv 1$ (mod $NN'$), we deduce
$\alpha'=\alpha$ and $\lambda'=\lambda$, which shows that $\lambda$
is obtained by the procedure described in Remark~\ref{complete_description}.
\end{proof}

\begin{proof}[Proof of Theorem~\ref{all_roots}]
We use the description of the roots of $b_{\fa}$ from
Remark~\ref{ingred}.
We can find $A\subseteq\{1,\ldots,r\}$, $B\subseteq\{1,\ldots,n\}$ and
$(\alpha_j)_{j\in A}$ in $\ZZ_{\geq 0}^{|A|}$ and 
$(\beta_i)_{i\in B}$ in $\ZZ_{<0}^{|B|}$ 
such that 
$$\Gamma:=\{w=(w_j)\in\QQ^r\mid w_j=\alpha_j\,{\rm for}\,j\,
{\rm in}\, A\,{\rm and}\,\ell_i(w)=\beta_i\,
{\rm for}\,i\,{\rm in}\, B\}$$
is nonempty, and if $w\in\Gamma$, then $\lambda=\sum_{j=1}^rw_j$. Moreover, 
we may assume by Lemma~\ref{condition} that the following condition holds:
if $u=(u_j)$ in $\ZZ^r$ is such that $u_j\geq-\alpha_j$ for all $j$ in $A$
and $\ell_i(u)\leq-\beta_i-1$ for all $i$ in $B$,
then $\sum_{j=1}^ru_j\leq 0$.

Let us fix $w$ in $\Gamma$. We deduce from our condition on $(\alpha_j)$ and 
$(\beta_i)$ that
\begin{equation}\label{lambda}
\lambda=\max_v\sum_{j=1}^rv_j,
\end{equation}
the maximum being over 
those $v$ in $\QQ^r$ such that $v_j\geq 0$ for $j$ in $A$, $\ell_i(v)\leq -1$
for $i$ in $B$
and $v_j-w_j\in\ZZ$ for all $j$ (the maximum is achieved for $v=w$).

We choose now $u=(u_j)$ in $\QQ_+^r$ such that 
$u_j=0$ if and only if $j$ is in $A$,
and such that $u+w$ is in $\ZZ^r$. Moreover, we may choose $u$ such that
$\ell_i(u)=\sum_{j=1}^ra_{i,j}u_j>0$ 
for all $i$ in $B$. Indeed, if this is not the case, then
there is $i$ in $B$ such that $a_{i,j}=0$ whenever $j$ is not in $A$.
This contradicts the fact that $\ell_i(w)<0$. Note also that 
since $u+w$ lies in $\ZZ^r$, $\ell_i(u)$ is an integer for every $i$ in $B$. 

Let $N'$ be a positive integer such that $N'u$ is in $\ZZ_{\geq 0}^r$.
If $q=mN'+1$ for $m\geq 1$, then $qu+w$  is in $\ZZ^r$.
With the notation introduced for the proof of Theorem~\ref{thm1},
we claim that if $q$ as above is large enough, then we have
\begin{equation}\label{claim}
\tau(B;(q\ell_i(u)-1)_{i\in B})=q\sum_{j=1}^ru_j+\lambda.
\end{equation}
This implies the assertion of the theorem:
take $J=(X_i^{\ell_i(u)}\vert i\in B)$ and $\alpha=\sum_{j=1}^ru_j$.

In order to prove the claim, suppose that $v\in\ZZ_{\geq 0}^r$ is such that
$\ell_i(v)\leq q\ell_i(u)-1$ for all $i$ in $B$. For $j$ in $A$ we have
$v_j-qu_j=v_j\geq 0$, and for $i$ in $B$ we have 
$\ell_i(v-qu)\leq -1$. As all $v_j-qu_j-w_j$ are integers, we deduce from
(\ref{lambda}) that $\sum_{j=1}^rv_j\leq q\sum_{j=1}^ru_j+\lambda$.

On the other hand, if $q\gg 0$ then $qu_j+w_j\geq 0$ for all $j$.
Note that $\ell_i(qu+w)\leq q\ell_i(u)-1$ for $i$ in $B$
and $\sum_{j=1}^r(qu_j+w_j)=q\sum_ju_j+\lambda$.
This completes the proof of the claim, and hence that of the theorem.
\end{proof}

We can prove now the description of the classes mod ${\mathbb Z}$
of the roots of the Bernstein-Sato polynomial.

\begin{proof}[Proof of Corollary~\ref{mod_Z}]
We use Theorems~\ref{thm1} and \ref{all_roots}. 
Recall that we have seen in Remark~\ref{rem_for_cor}
that the classes in $\QQ/\ZZ$ of the roots of $b_{\fa}$
are equal to 
\begin{equation}\label{eq100}
\{-L_{\sigma}(w)+\ZZ\mid\sigma\,{\rm maximal}\, {\rm cone}\, {\rm in}\,\Delta,
w\in\sigma\cap (\ZZ_{>0})^n\}.
\end{equation}
Therefore in order to prove the theorem it is enough to show that for every
maximal cone $\sigma$ in $\Delta$, the set of classes
$\{-L_{\sigma}(w)+\ZZ\mid w\in\sigma\cap (\ZZ_{>0})^n\}$ is the subgroup of
$\QQ/\ZZ$ generated by $\frac{1}{m_{\sigma}}$, where if $Q$ is the facet
of $P_{\fa}$ corresponding to $\sigma$, 
we put $m_{\sigma}$
for $m_{Q}$. 

Since $m_{\sigma}L_{\sigma}$ has integer coefficients, for every
$w$ in $\sigma\cap (\ZZ_{>0})^n$ the class of $-L_{\sigma}(w)$ lies in
the subgroup generated by $1/m_{\sigma}$. On the other hand, if
$\tau\subseteq\sigma$ is a convex cone generated by a basis 
$e'_1,\ldots,e'_n$ for $\ZZ^n$, then $m_{\sigma}$ is the smallest positive
integer such that all $m_{\sigma}L_{\sigma}(e'_i)$ are integers. 
In this case it follows easily that there is $w$ in the interior of
$\tau$ such that $\frac{1}{m_{\sigma}}+L_{\sigma}(w)$ is an integer.
By taking suitable multiples of $w$ we see that the subgroup 
generated by the class of $\frac{1}{m_{\sigma}}$ is contained 
in (\ref{eq100}),
which completes the proof.
\end{proof}

\begin{example}
As we have already mentioned, it is a general fact that the largest
root of the Bernstein-Sato
polynomial of $\fa$ is $-\lct(\fa)$, where $\lct(\fa)$ is the log canonical
threshold of $\fa$ (see \cite{BMS}). 
Let us prove this for monomial ideals using the
above results.

It follows from \cite{Ho} that $\lct(\fa)$ is the
largest positive real number $c$ such that 
$e=(1,\ldots,1)$ lies in $cP_{\fa}$, where $P_{\fa}$ is the
Newton polyhedron of $\fa$.
It follows from Lemma~{\ref{lemma1}} that, with our notation,
 $\lct(\fa)=\tau_{\QQ}(e)$. 

Let $\sigma_0$ be a maximal cone
in the fan $\Delta$
such that $e$ lies in $\sigma_0$. With the notation in 
Remark~\ref{complete_description}, we get the root $-\lct(\fa)$
as $L_{\sigma_0}(-e)-A^{\sigma_0}_0$ (corresponding to $c=(0,\ldots,0)$). 

On the other hand, if $\lambda$ is another root of $b_{\fa}$,
then by Theorem~\ref{all_roots} and Remark~\ref{complete_description}
we have $\lambda=L_{\widetilde{\sigma}}(-e)-A^{{\sigma}}_c$
for some $\sigma$, $\widetilde{\sigma}$ and some $c$. 
Since $A^{{\sigma}}_c\geq 0$ and since concavity of 
$\tau_{\QQ}$ gives $L_{\widetilde{\sigma}}(e)
\geq \tau_{\QQ}(e)=L_{\sigma_0}(e)$,
we get $\lambda\leq-\lct(\fa)$.
\end{example}

\bigskip

We end this section with 
a description of the $F$-thresholds of monomial ideals
(see \S 2 for the definition). 
It follows from Proposition~\ref{prop_A1} in the Appendix
that the set of all $F$-thresholds of the monomial ideal $\fa$
(computed for the reduction mod $p$ of $\fa$, where $p\gg 0$
is a prime) is equal to the set of $F$-thresholds $c^J(\fa_p)$
with respect to monomial ideals $J$. Moreover, we have seen
that it is enough to consider the case when
$J=(X_1^{b_1},\ldots,X_n^{b_n})$,
where $b_i$ are positive integers.  
With the notation in the proof
of Theorem~\ref{thm1}, we have shown that for this $J$ we have
$$\pt^{J}(\fa_p)=\tau_{\QQ}(b)$$
for every prime $p$ large enough. We will denote this number simply by 
$\pt^J(\fa)$, as it does not depend on $p$. 

We make the connection with the multiplier ideals of $\fa$. Recall that by the
description in \cite{Ho}, the multiplier ideal $\cI(\fa^{\alpha})$
of $\fa$ with exponent $\alpha$ can be described as follows
$$\cI(\fa^{\alpha})=(X^w\vert u+e\in {\rm Int}(\alpha P_{\fa})).$$

The jumping coefficients of the multiplier ideals defined in \cite{ELSV}
are those $\alpha>0$ such that $\cI(\fa^{\alpha})$ is strictly
contained in $\cI(\fa^{\alpha-\epsilon})$ for every $\epsilon>0$.
It follows from the above description of the multiplier ideals of $\fa$,
that every $b=(b_i)$ with $b_i$ positive integers gives a jumping coefficient
$\alpha$ characterized by 
the fact that $b$
lies in the boundary of $\alpha P_{\fa}$. Moreover, every 
jumping coefficient arises in this way.
The connection with the $F$-thresholds is given by

\begin{proposition}\label{desc}
For every prime $p$,
the jumping coefficient corresponding to $b=(b_i)$ as above is equal to
$\pt^J(\fa_p)$, where $J=(X_1^{b_1},\ldots,X_n^{b_n})$.
\end{proposition}

\begin{proof}
The assertion follows from Lemma~{\ref{lemma1}}.
\end{proof}

\begin{remark}
It follows from the above discussion that the $F$-thresholds in the monomial case
are easy to describe. Alternatively, this can be seen as follows:
it is shown in \cite{MTW} that 
in general, the $F$-thresholds (for various $J$) are the jumping
coefficients for the test ideals introduced in \cite{HY} as an analogue
of multiplier ideals. On the other hand, it is proved in \cite{HY} that
for monomial ideals the test ideals coincide with the multiplier ideals.
This gives a different approach to the description in Proposition~\ref{desc}. 
\end{remark}

\begin{remark}
It is clear that for monomial ideals all $F$-thresholds are rational numbers.
In fact, more is true. If $p$ is a prime, then
the series $\sum_e\nu^J_{\fa}(p^e)t^e$ is rational. This follows from the fact
that the function $e\to \nu^J_{\fa}(p^{e+1})-p\nu^J_{\fa}(p^e)$ 
is eventually periodic, which
is a consequence of  Theorem~\ref{thm1}. 
\end{remark}

\section{Examples}

In this section we give some examples
to illustrate how to use our approach
to give roots of the Bernstein-Sato polynomial.
We use freely the notation introduced for the proof of Theorem~\ref{thm1}.
In all these examples we do not describe the complete
picture given in Remark~\ref{complete_description}, but 
we give enough information
to recover all the roots.
For more complicated examples, based on a more explicit
combinatorial description of the roots of the Bernstein-Sato polynomial,
we refer to \cite{BMS1}.

\begin{example}\label{ex3}
Let $\fa=(\prod_{j\neq i}X_j\vert 1\leq i\leq 4)$ in $\ZZ[X_1,\ldots,X_4]$,
so we have
$\ell_i(s)=\sum_{j\neq i}s_j$ for $1\leq i\leq 4$.
If $\alpha\in\NN^4$ is such that $\ell_i(\alpha)\leq b_i$ for all $i$,
summing these inequalities we get
$$\sum_i\alpha_i\leq\lfloor (\sum_ib_i)/3\rfloor,$$
where we use the notation $\lfloor x\rfloor$ for the largest integer
$\leq x$.

Consider first the cone 
$$\sigma=\{w\vert w_i\geq 0, \sum_jw_j\geq 3w_i\,{\rm for}\,{\rm all}\,i\}.$$
If $b\in\sigma\cap\NN^4$, then $\tau(b)=\lfloor (\sum_ib_i)/3\rfloor$.
Indeed, suppose for example that $\sum_ib_i\equiv 2$ (mod $3$). 
If we take $\alpha_i=(\sum_jb_j-2)/3-(b_i-1)$ for $i=1,2$ and
$\alpha_i=(\sum_jb_j-2)/3-b_i$ for $i=3,4$, then $\alpha\in\NN^4$
and $\sum_i\alpha_i=(\sum_ib_i-2)/3$. The other two cases are similar.
We get $L_i$ for $i=0,1,2$ with $L_i(b)=(\sum_jb_j-i)/3$, such that
whenever $b\in\sigma\cap\NN^4$, with $\sum_ib_i
\equiv i$ (mod $3$) we have $\tau(b)=L_i(b)$.
All $L_i(-1,\ldots,-1)$ are roots of $b_{\fa}$, which gives
the roots $-\frac{4}{3}$,$-\frac{5}{3}$ and $-2$.

Suppose now that $b\in\NN^4$ is such that, for example,
$2b_4>b_1+b_2+b_3$.
If $\ell_i(\alpha)\leq b_i$ for all $i$, adding these inequalities for 
$1\leq i\leq 3$ implies
$$\sum_{i=1}^4\alpha_i\leq\lfloor(\sum_{i=1}^3b_i)/2\rfloor.$$
If we assume, in addition, that $b_1+b_2\geq b_3$, $b_1+b_3\geq b_2$
and $b_2+b_3\geq b_1$, then
this maximum can be achieved. Indeed, if $\sum_{i=1}^3b_i$ is even,
take $\alpha_j=(\sum_{i=1}^3b_i)/2-b_j$ for $1\leq j\leq 3$ and $\alpha_4=0$.
If $\sum_{i=1}^3b_i$ is odd, take $\alpha_1=(\sum_{i=1}^3b_i-1)/2-(b_1-1)$,
$\alpha_j=(\sum_{i=1}^3b_i-1)/2-b_j$ for $j=2,3$, and $\alpha_4=0$.

If we take $L'_i$ for $i=0,1$ given by $L'_i(b)=(\sum_{j=1}^3b_j-i)/2$,
we see that if $b$ is as above and $\sum_{j=1}^3b_j\equiv i$
(mod $2$), we have $\tau(b)=L'_i(b)$. Therefore both $L'_i(-1,\ldots,-1)$
are roots of $b_{\fa}$, which gives the roots
$-\frac{3}{2}$ and $-2$. 

Therefore we have obtained the roots $-\frac{3}{2}$, $-\frac{4}{3}$,
$-\frac{5}{3}$ and $-2$.  Note that, in fact, by \cite{BMS} (4.5)
$b_{\fa}=(s+\frac{3}{2})(s+\frac{4}{3})(s+\frac{5}{3})(s+2)^3$.
\end{example}

\begin{example}\label{ex2}
Let $\fa=(X^2YZ,XY^2Z,XYZ^2)$, so
$\ell_1(s)=2s_1+s_2+s_3$, $\ell_2(s)=s_1+2s_2+s_3$ and
$\ell_3(s)=s_1+s_2+2s_3$. 

If $\alpha\in\NN^3$ is such that
$\ell_i(\alpha)\leq b_i$, by summing these relations we get
$\sum_i\alpha_i\leq\lfloor\frac{1}{4}\sum_ib_i\rfloor$.
This maximum is achieved if $b$ lies in the cone
$$\sigma=\{w\vert w_i\geq 0, 
w_i\geq\frac{\sum_jw_j}{4}\,{\rm for}\,{\rm all}\,i\}.$$
Indeed, suppose for example that $\sum_jb_j\equiv 1$ (mod $4$). 
If $\alpha_1=b_1-\frac{\sum_ib_i-1}{4}$, $\alpha_2=b_2-\frac{\sum_ib_i-1}{4}$,
and $\alpha_3=b_3-\frac{\sum_ib_i+3}{4}$, then $\alpha\in\NN^3$ and $\sum_i\alpha_i
=(b_1+b_2+b_3-1)/4$. The other three cases are similar.

We get $L_i(b)=(\sum_jb_j-i)/4$ for $i=0,1,2,3$ such that
if $b\in\NN^3$ is in  $\sigma$ and $\sum_jb_j\equiv i$ (mod $4$),
then $\tau(b)=L_i(b)$.
We get roots  of the Bernstein-Sato polynomial of $\fa$ given
by $L_i(-1,-1,-1)$ for $i=0,1,2,3$, i.e. $-\frac{3}{4}$, $-1$, $-\frac{5}{4}$
and $-\frac{6}{4}$. Note that by \cite{BMS} (4.5), 
$b_{\fa}(s)=(s+\frac{3}{4})(s+\frac{5}{4})(s+\frac{6}{4})(s+1)^3$.
\end{example}

\begin{example}\label{ex1}
Let $\fa$ be the ideal generated by $X_iX_j$ for all $1\leq i<j\leq n$,
with $n\geq 3$. We have
$$\ell_i(s)=\sum_{j<i}s_{j,i}+\sum_{j>i}s_{i,j}$$
for $s=(s_{i,j})_{i<j}$.
Note that if $\ell_i(s)\leq b_i$ for all $i$, by taking the sum we get
$$\sum_{i<j}s_{i,j}\leq\lfloor (\sum_ib_i)/2\rfloor.$$

It is easy to give, as above, a full-dimensional cone $\sigma$ on which this
maximum is achieved. This gives $L_1$ and $L_2$ with
$L_i(b)=(\sum_jb_j-i)/2$ such that if $b\in\NN^n$
lies in $\sigma$ and $\sum_jb_j\equiv i$
(mod $2$), then $\tau_i(b)=L_i(b)$. 
Therefore $L_0(-1,\ldots,-1)=-\frac{n}{2}$ and $L_1(-1,\ldots,-1)=-\frac{n+1}{2}$
are roots of $b_{\fa}$.

Suppose now that $b\in\NN^n$ is such that $b_1>\sum_{j=2}^nb_j$.
If $\ell_i(\alpha)\leq b_i$ for all $i$, by adding these inequalities for
$i\geq 2$, we deduce $\sum_{i<j}\alpha_{i,j}\leq\sum_{j\geq 2}b_j$.
Moreover, this maximum can be achieved: let $\alpha_{1,i}=b_i$ for 
$i>1$ and all other $\alpha_{i,j}=0$.
If $L(b)=\sum_{j\geq 2}b_j$ we see that if $b$ is as above, then
$\tau(b)=L(b)$. Hence $L(-1,\ldots,-1)=-(n-1)$ is also a root of
$b_{\ba}$. Note that by \cite{BMS} (4.5), the Bernstein-Sato
polynomial of $\fa$ is $b_{\fa}(s)=(s+\frac{n}{2})
(s+\frac{n+1}{2})(s+n-1)$.
\end{example}

\section{Appendix}

We show that if $\fa$ is a monomial ideal, then in order to compute
the functions $\nu_{\fa}^J(p^e)$ for $p\gg 0$, it is enough to consider
the case when $J$ is a monomial ideal. 

\begin{proposition}\label{prop_A1}
Let $\fa$ be a nonzero ideal generated by
monomials in $\ZZ[X_1,\ldots,X_n]$ and let 
$J\subseteq (X_1,\ldots,X_n)$
be an ideal such that $\fa$ is contained in the radical of $J$. If
we put
$$\widetilde{J}:=(X^u\mid h X^u\in J\,{\rm for}\,{\rm some}\, 
h\in\ZZ[X_1,\ldots,X_n]\smallsetminus (X_1,\ldots,X_n))$$
and if $p\gg 0$, then $\nu^J_{\fa}(p^e)=\nu_{\fa}^{\widetilde{J}}(p^e)$
for every $e\geq 1$.
\end{proposition}

For the proof of Proposition~\ref{prop_A1} we will need the following lemma.

\begin{lemma}\label{lem_A2}
With the notation in Proposition~\ref{prop_A1},
for $p\gg 0$ we have
$$\widetilde{J}_p=(X^u\mid X^u\in J_p)$$
as ideals in ${\mathbb F}_p[X_1,\ldots,X_n]_{(X_1,\ldots,X_n)}$.
\end{lemma}

\begin{proof}
Given $X^u$ in $\widetilde{J}$ and $h\in\ZZ[X_1,\ldots,X_n]$
such that $h(0)\neq 0$ and $hX^u\in J$, then for $p$ not dividing
$h(0)$ we see that the image of $X^u$ 
in ${\mathbb F}_p[X_1,\ldots,X_n]_{(X_1,\ldots,X_n)}$
lies in $J_p$. 
Since $\widetilde{J}$ is finitely generated, for $p\gg 0$ we deduce the
inclusion $\subseteq$ in the statement.

We prove now the reverse inclusion.
It is clear that if $a=(a_i)$ and $b=(b_i)$ in $\NN^n$ are such that
$a_i\leq b_i$ for every $i$, then $(J\colon X^a)\subseteq
(J\colon X^b)$. For a subset $A$ of $\{1,\ldots,n\}$ and $w=(w_i)\in\NN^n$,
we denote by $w+\NN^{A}$ the set 
$$\{u=(u_i)\in\NN^n\mid u_i\geq w_i\,{\rm for}
\,{i\in A}\,{\rm and}\,u_i=w_i\,{\rm for}\,i\not\in A\}.$$
 
Using the fact that $\ZZ[X_1,\ldots,X_n]$ is Noetherian, 
we deduce that there is a decomposition
$$\NN^n=\bigsqcup_{j=1}^m(w^{(j)}+\NN^{A_j}),$$
for some $w^{(j)}$ in $\NN^n$ and $A_j\subseteq\{1,\ldots,n\}$ such that
for every $w$ in $w^{(j)}+\NN^{A_j}$ the ideal $(J\colon X^w)$ is equal 
to a fixed ideal $I_j$.
In particular, we deduce that for every $i$ in $A_j$, we have
$(I_j\colon X_i)=I_j$. 

Given a prime $p$,
for every ideal $I$ in $\ZZ[X_1,\ldots,X_n]$ we denote by 
$\overline{I}$ the reduction of $I$ in ${\mathbb F}_p[X_1,\ldots,X_n]$.
There is $p_0$ such that for every prime $p\geq p_0$ the following hold:
\begin{equation}\label{eq5}
(\overline{I_j}\colon X_i)=\overline{I_j}\,\,{\rm in}
\,{\mathbb F}_p[X_1,\ldots,X_n]\,\,{\rm for}\,{\rm all}\,
j\leq m\,{\rm and}\,{\rm all}\,i\in A_j,
\end{equation}
\begin{equation}\label{eq6}
(\overline{J}\colon X^{w^{(j)}})=\overline{I_j}\,\,{\rm in}\,
{\mathbb F}_p[X_1,\ldots,X_n]\,\,{\rm for}\,{\rm all}\,j\leq m.
\end{equation}
In order to show (\ref{eq5}), use the exact sequence
\begin{equation}\label{eq10}
0\to (I_j\colon X_i)/I_j\to \ZZ[X_1,\ldots,X_n]/I_j
\overset{h}\to \ZZ[X_1,\ldots,X_n]/I_j,
\end{equation}
where $h$ is multiplication by $X_i$.
Note that given an arbitrary finitely generated $\ZZ[X_1,\ldots,X_n]$-module
$M$, there is a positive integer $\ell$ such that $M[1/\ell]$
is flat over $\ZZ[1/\ell]$.
Indeed, the submodule
 $N\subseteq M$ of elements annihilated by some positive integer
 is finitely generated over $\ZZ[X_1,\ldots,X_n]$. Therefore we can find $\ell$
such that $\ell N=0$, and $M[1/\ell]$ is torsion-free, hence flat over
$\ZZ[1/\ell]$. Applying this observation for the image and for the cokernel
of $h$, we see that if $p\gg 0$, then
the sequence obtained 
 from (\ref{eq10})
by tensoring with ${\mathbb F}_p$ is again exact, which gives (\ref{eq5}).
The proof of (\ref{eq6}) is similar.

Consider now a prime $p\geq p_0$ and fix $w\in\NN^n$. Let $j$ be such that
$w$ is in $w^{(j)}+\NN^{A_j}$. Suppose that $h$ is in 
${\mathbb F}_p[X_1,\ldots,X_n]$ such that $h(0)\neq 0$ and
$hX^w$ is in $\overline{J}$. Using (\ref{eq6}) we deduce that
$hX^{w-w_j}$ lies in $\overline{I_j}$. Moreover, (\ref{eq5})
implies that $h$ is in $\overline{I_j}$. If $g$ is an element in $I_j$
whose class is $h$, then $g(0)\neq 0$ and $gX^{w_j}$ is in $J$. In particular,
$gX^w$ is in $J$, so $X^w\in \widetilde{J}$,
which completes the proof of the lemma.
\end{proof}

\begin{proof}[Proof of Proposition~\ref{prop_A1}]
We consider $p\gg 0$ so the assertion in Lemma~\ref{lem_A2} applies. 
In particular, we have
$\widetilde{J}_p\subseteq J_p$, so $\nu^{\widetilde{J}}_{\fa}(p^e)
\geq\nu^{J}_{\fa}(p^e)$. In order to show that we have equality,
since $\fa$ is monomial it is enough to prove that if $X^u$ is in
$J_p^{[p^e]}$, then $X^u$ is in $\widetilde{J}_p^{[p^e]}$. 
We do induction on $e\geq 0$, the case $e=0$ being a consequence
of Lemma~\ref{lem_A2}.

If $e\geq 1$, we write $u=pv+w$ with $v$ and $w$ in $\NN^n$ and
$0\leq w_i<p$ for every $i$. We deduce 
\begin{equation}\label{e_case}
X^w\in (J_p^{[p^e]}\colon X^{pv})=(J_p^{[p^{e-1}]}\colon X^v)^{[p]},
\end{equation}
where the equality follows from the fact that ${\mathbb F}_p[X_1,\ldots,X_n]$
is flat over ${\mathbb F}_p[X_1^p,\ldots,X_n^p]$.
The assumption on $w$ implies that $w_i=0$ for every $i$ and 
$(J_p^{[p^{e-1}]}\colon X^v)$ is the unit ideal. 
We deduce that
 $X^v$ lies in $J_p^{[p^{e-1}]}$, hence in $\widetilde{J}_p^{[p^{e-1}]}$
by the induction hypothesis. Therefore $X^u$ is in $J_p^{[p^e]}$,
which completes the proof of the proposition.
\end{proof}

\providecommand{\bysame}{\leavevmode \hbox \o3em
{\hrulefill}\thinspace}


\begin{thebibliography}{MTW}

\bibitem[Bj]{Bj}
J.-E. Bj\"{o}rk, \emph{Rings of differential operators}, 
Amsterdam, North-Holland, 1979.

\bibitem[BMS1]{BMS1}
N.~Budur, M.~Musta\c{t}\v{a} and M.~Saito,
Combinatorial description of the 
roots of the Bernstein-Sato polynomials for monomial ideals,
preprint 2005. 

\bibitem[BMS2]{BMS}
N.~Budur, M.~Musta\c{t}\v{a} and M.~Saito,
Bernstein-Sato polynomials of arbitrary varieties,
math.AG/0408408.

\bibitem[ELSV]{ELSV}
L.~Ein, R.~Lazarsfeld, K.~E.~Smith and D.~Varolin,
Jumping coefficients of multiplier ideals, Duke Math. J. \textbf{123} (2004),
469--506.


\bibitem[HY]{HY}
N.~Hara and K.-i.~Yoshida, A generalization of tight closure and multiplier 
ideals, Trans. Amer. Math. Soc. \textbf{355} (2003), 3143--3174.


\bibitem[Ho1]{Ho}
J.~Howald, Multiplier ideals of monomial ideals,
Trans. Amer. Math. Soc. \textbf{353} (2001), 2665--2671.

\bibitem[Ka]{kashiwara}
M.~Kashiwara, \emph{$D$-modules and microlocal calculus}, 
translated from the 2000 Japanese original by Mutsumi Saito,
Translations of Mathematical
   Monographs \textbf{217}, Iwanami Series in Modern Mathematics,
American Mathematical Society, Providence, RI, 2003.


\bibitem[Ko]{kollar}
J.~Koll\'{a}r,
Singularities of pairs, in \emph{Algebraic geometry,
Santa Cruz 1995}, volume \textbf{62} of Proc. Symp. Pure Math Amer.
Math. Soc. 1997, 221--286.

\bibitem[La]{lazarsfeld}
R.~Lazarsfeld, \emph{Positivity in algebraic geometry II}, 
Ergebnisse der Mathematik und ihrer Grenzgebiete. 3. Folge, 
A series of Modern Surveys in Mathematics, Vol. \textbf{49},
Springer-Verlag, Berlin, 2004.

\bibitem[MTW]{MTW}
M.~Musta\c{t}\v{a}, S.~Takagi and K.-i.~Watanabe,
$F$-thresholds and Bernstein-Sato polynomials, 
to appear in A.~Laptev (ed.), \emph{European congress of mathematics
(ECM), Stockolm, Sweden, June 27--July 2, 2004},
Z\"{u}rich, European Mathematical Society, 2005.

\bibitem[Zi]{Zi}
G.~Ziegler,
 Lectures on polytopes,
Graduate Texts in Mathematics \textbf{152},
Springer-Verlag, New York, 1995.


\end{thebibliography}
\end{document}